\documentclass[12pt,twoside]{article}

\usepackage{amsmath,amssymb,amsthm} 

\newcommand{\p}{\mathfrak{p}}
\newcommand{\Prim}{\mathrm{Prim}}
\newcommand{\Conj}{\mathrm{Conj}}
\newcommand{\vol}{\mathrm{vol}}

\newcommand{\Id}{\mathrm{Id}}

\newcommand{\Gal}{\mathrm{Gal}}

\newcommand{\ord}{\mathrm{ord}}
\newcommand{\as}{\quad\text{as}\quad}
\newcommand{\tinf}{\to\infty}
\newcommand{\disp}{\displaystyle}
\newcommand{\bsla}{\backslash}

\newcommand{\bC}{\mathbb{C}}
\newcommand{\bR}{\mathbb{R}}
\newcommand{\bQ}{\mathbb{Q}}
\newcommand{\bZ}{\mathbb{Z}}
\newcommand{\bH}{\mathbb{H}}

\newcommand{\noi}{\noindent}
\renewcommand{\Re}{\mathrm{Re}}
\renewcommand{\Im}{\mathrm{Im}}
\newcommand{\vphi}{\varphi}
\newcommand{\divset}{\hspace{3pt}|\hspace{3pt}}

\newtheorem{thm}{Theorem}[section]

\newtheorem{lem}[thm]{Lemma}
\newtheorem{rem}[thm]{Remark}

\newtheorem{claim}{Claim}[section]

\numberwithin{equation}{section}

\title{Partial zeta functions}
\author{Yasufumi Hashimoto \thanks{Supported by JSPS Research Fellowships for Young Scientists}}
\date{}

\begin{document}
\markboth
{Y. Hashimoto}
{Partial zeta functions}
\pagestyle{myheadings}

\maketitle
\begin{abstract}
In this paper, we study analytic properties of zeta functions defined by partial Euler products.
\end{abstract}

\renewcommand{\thefootnote}{}

\section{Introduction}
In the previous paper \cite{HW}, we have studied the splitting densities of 
prime geodesics of negatively curved locally symmetric Riemannian manifolds $X$ 
in a finite cover $\tilde{X}$ of $X$ as an extension of the prime geodesic theorem.
Especially, when the fundamental group of $X$ is $SL_2(\bZ)$ and 
that of $\tilde{X}$ is a congruence subgroup of $SL_2(\bZ)$, 
we have explicitely determined the type of splitting for each geodesic, 
and calculate the splitting densities for every types. 
Applying the results in \cite{HW} to the formula in \cite{VZ} 
about the relation between the Selberg zeta functions for $X$ and $\tilde{X}$, 
we can obtain an expression of the Selberg zeta function for $\tilde{X}$ 
as a product over prime geodesic of $X$. 
By taking the quotients of such expressions of two Selberg zeta functions, 
we find a formula of the zeta function given by the Euler product over prime geodesics 
with a certain type of splitting. 
From this formula, we have obtained the analytic continuation to $\{\Re{s}>0\}$ 
of such a partial Selberg zeta function (see \cite{HW}).

On the other hand, in \cite{Ku}, it was studied 
the zeta functions defined by the Euler products over prime numbers $p$ 
satisfying $(d/p)=1$ (or $=-1$) for a fixed square free integer $d$.
In fact, it was shown that these zeta functions are analytically continued to $\{\Re{s}>0\}$ 
and have natural boundaries on $\Re{s}=0$.

The aim of the present paper is to generalize zeta functions defined by partial Euler products 
and to study their analytic properties. 
We first get in Theorem \ref{thm1} the analytic continuations in $\{\Re{s}>0\}$ 
by extending the idea used in \cite{Ku}. 
Furthermore, in Theorem \ref{thm2}, we state the sufficient condition of the distributions 
of non-trivial zeros 
for the partial zeta functions having natural boundaries on $\Re{s}=0$. 
As examples, we treat the cases of the Dedekind zeta functions, the Selberg zeta functions 
and the Ihara zeta functions of graphs.  
Actually we show that the partial zeta functions of the Dedekind zeta functions for abelian extensions 
of $\bQ$, 
the Selberg zeta functions for congruence subgroups of $SL_2(\bZ)$ and the Ihara zeta functions 
for finite regular graphs have natural boundaries on $\Re{s}=0$.

\section{Notations and main results}
Let $P$ be an infinite countable set and $N:P\to\bR_{>1}$ a map such that $\sum_{p\in P}N(p)^{-d}<\infty$
for some $d>0$. 
Put $d_P:=\inf\{d>0\divset\sum_{p\in P}N(p)^{-d}<\infty\}$ and assume that $d_P>0$. 
For convenience, we normalize $N$ by $d_P=1$.
We define the zeta function of $P$ by 
\begin{align*}
\zeta_P(s):=\prod_{p\in P}(1-N(p)^{-s})^{-1}\quad \Re{s}>1,
\end{align*}
and assume that   
(i) $\zeta_P(s)$ is non-zero holomorphic in $\{\Re{s}> 1\}$ and has a simple pole at $s=1$, and 
(ii) $\zeta_P(s)$ has an analytic continuation to the whole complex plane $\bC$ as a meromorphic function. 

Let $G$ be a finite group and $\hat{G}$ the set 
of the finite dimensional irreducible unitary representations of $G$.
For a map $\vphi:P\to\Conj(G)$ and $\rho\in\hat{G}$, we define the $L$-functions by
\begin{align*}
L_P^{(G)}(s,\rho)=L_P(s,\rho):=\prod_{p\in P}\det{(1-\rho(\vphi(p))N(p)^{-s})}^{-1}\quad \Re{s}>1,
\end{align*}
and assume that, if $\rho\not=1$, then 
$L_P(s,\rho)$ satisfies that (i') $L_P(s,\rho)$ is non-zero holomorphic in $\{\Re{s}\geq 1\}$ 
and satisfies the same condition (ii) for $\zeta_P(s)$.

Put $P_n(G)=P_n:=\{p\in P\divset\ord_G(\vphi(p))=n\}$, where $\ord_G(\vphi(p))$ is the order 
of $\vphi(p)$ in $G$.
In the present paper, we study analytic properties of the following zeta function.
\begin{align*}
\zeta_{P_n}(s):=\prod_{p\in P_n}(1-N(p)^{-s})^{-1}\quad \Re{s}>1.
\end{align*}
For simplicity we treat the case where $G$ is a cyclic group of prime order $q$. 
First we get the following result.
\begin{thm}\label{thm1}
The function $\zeta_{P_q}(s)$ satisfies the following functional equations.
\begin{align}
\frac{\big(\zeta_{P_q}(s)\big)^q}{\zeta_{P_q}(qs)}=&\frac{\big(\zeta_P(s)\big)^q}{Z_P(s)}\label{Pq},
\end{align}
where $Z_P^{(G)}(s)=Z_{P}(s):=\prod_{\rho\in\hat{G}}L_P(s,\rho)^{\dim\rho}$.
Furthermore, for any $r\geq1$, the function $\big(\zeta_{P_q}(s)\big)^{q^r}$ is analytically continued
to $\{\Re{s}>1/q^r\}$ as a meromorphic function 
and has infinitely many singular points near $s=0$.
\end{thm}

\begin{proof}
By elementary calculations, we obtain 
\begin{align*}
Z_P(s)=&\prod_{n\mid\#G}\prod_{\rho\in\hat{G}}\prod_{p\in P_n}
\det(1-\rho(\vphi(p))N(p)^{-s})^{-\dim\rho}\\
=&\prod_{n\mid\#G}\prod_{p\in P_n}(1-N(p)^{-ns})^{-\#G/n}\\
=&\prod_{n\mid\#G}\big(\zeta_{P_n}(ns)\big)^{\#G/n}.
\end{align*}
Then, when $G$ is a cyclic group of order $q$, we have
\begin{align*}
Z_{P}(s)=&\big(\zeta_{P_1}(s)\big)^q\zeta_{P_q}(qs).
\end{align*}
Since $\zeta_P(s)=\zeta_{P_1}(s)\zeta_{P_q}(s)$, the equation \eqref{Pq} follows immediately. 

Now, for convenience, we rewrite the equation \eqref{Pq} as follows. 
\begin{align}
\frac{f(s)^q}{f(qs)}=g(s).\label{simplePq}
\end{align}
From \eqref{simplePq}, we recursively obtain the following formula.
\begin{align}
f(s)^{q^r}=&f(q^{r}s)\prod_{i=0}^{r-1}g(q^{i}s)^{q^{r-i-1}}.\label{analytic}
\end{align}
Since $g(s)$ is meromorphic and $f(q^{r}s)$ is non-zero holomorphic in $\{\Re{s}>1/q^r\}$, 
we see that $f(s)^{q^r}$ is meromorphic in $\{\Re{s}>1/q^r\}$. 

Since $g(s)=\zeta_{P}(s)^q/Z_{P}(s)$, the function $g(s)$ is non-zero holomorphic in $\{\Re{s}> 1\}$ 
and has a pole at $s=1$ of order $q-1$. 
Then, according to \eqref{analytic}, we see that $f(s)$ has branch points at $s=1/q^{i}$ for $i=0,1,\cdots$. 
This completes the proof of Theorem \ref{thm1}. 
\end{proof}

The theorem above gives the analytic continuation to $\{\Re{s}>0\}$ of $\zeta_{P_q}(s)$. 
Next we study analytic properties in $\{\Re{s}\leq0\}$ of $\zeta_{P_q}(s)$. 

Let $\Lambda$ be the set of singular points in $\{0<\Re{s}<1,\Im{s}>0\}$ of $g(s)$, 
and $m(\sigma)$ is the order of $\sigma\in\Lambda$ 
(when $\sigma$ is a pole, $m(\sigma)$ is a negative value). 
Denote by $\Omega:=\{q^{-k}\sigma\divset\sigma\in\Lambda,k\geq0\}$.
According to \eqref{analytic}, the set of singular points in $\{0<\Re{s}<1\}$ of $f(s)$ 
is a subset of $\Omega$.
For $\sigma\in\Lambda$, we denote by $[\sigma]$ the subset of $\Lambda$ 
which consists of elements $\sigma'=q^k\sigma\in\Lambda$ for some $k\in\bZ$.
We also denote by $M_q(\sigma):=\sum_{\sigma'\in[\sigma]}q^{-k}m(\sigma')$ 
and $\Lambda_q:=\{[\sigma]\subset\Lambda\divset M_q(\sigma)\neq0\}$. 
Under such notations, we obtain the following result.

\begin{thm}\label{thm2}
Number the elements of $\Lambda_q$ by $\sigma_0,\sigma_1,\sigma_2,\cdots$ such that 
$\sigma_i\neq\sigma_j$ for $i\neq j$ and $0<\beta_i\leq\beta_j$ for $i<j$, where $\beta_i:=\Im{\sigma_i}$. 
If $\beta_j\tinf$ and $(\beta_j)^{1/j}\to1$ as $j\tinf$, 
then the partial zeta function $\zeta_{P_q}(s)$ has a natural boundary on $\Re{s}=0$.
\end{thm}

\begin{proof}
According to \eqref{analytic}, we see that the singular points of $f(s)$ 
near the line $\Re{s}=0$ consist in 
\begin{align*}
\Omega_q:=\{q^{-k}\sigma\divset \sigma\in\Lambda_q,k\geq0\}.
\end{align*}
We now assume that $f(s)$ does not have a natural boundary on $\Re{s}=0$, 
namely there exist constants $T_1,T_2>0$ ($T_2>T_1$) 
such that $\beta_jq^{-k}<T_2$ or $\beta_jq^{-k}>T_1$ for any $j,k\geq0$. 
Put $j(T):=\{j\geq0\divset\beta_j\leq T\leq\beta_{j+1}\}$ for $T>0$, 
and $J(k):=j(T_1q^k)$ for a fixed $k\geq0$.
By the assumption, we have $q^{-k}\beta_{J(k)+1}>T_2$. 
It is easy to see that $j(T)$ and $J(k)$ are non-decreasing functions of $T$ and $k$ respectively.  

We now estimate $\beta_{J(k)+1}-\beta_{J(k)}$ for sufficiently large $k>0$. 
By the definition of $J(k)$ and the assumption, we have
\begin{align}
q^{-k}(\beta_{J(k)+1}-\beta_{J(k)})>T_2-T_1>0.\label{estimate1}
\end{align}
On the other hand, since $(\beta_j)^{1/j}\to1$, we have
\begin{align*}
\beta_{j+1}-\beta_j=o(\beta_j)\as j\tinf.
\end{align*}
Then we obtain 
\begin{align}
q^{-k}(\beta_{J(k)+1}-\beta_{J(k)})=q^{-k}o(\beta_{J(k)})<q^{-k}o(T_1q^{k})=o(1)\as k\tinf.\label{estimate2}
\end{align}
The estimates \eqref{estimate1} and \eqref{estimate2} contradict to each other. 
Then the assumption is false and, therefore, Theorem \ref{thm2} holds.
\end{proof}

\begin{rem}
It is not difficult to obtain results such as Theorem \ref{thm1} and \ref{thm2} 
for the case that $G$ is not necessarily a cyclic group of prime order. 
For example, when $\#G=q_1q_2$ for distinct primes $q_1$ and $q_2$, we can obtain 
\begin{align}
\frac{\zeta_{P_{q_1q_2}}(q_1q_2s)\zeta_{P_{q_1q_2}}(s)^{q_1q_2}}{\zeta_{P_{q_1q_2}}(q_1s)^{q_2}
\zeta_{P_{q_1q_2}}(q_2s)^{q_1}}
=\frac{Z_P^{(G)}(s)\zeta_{P}(s)^{q_1q_2}}{Z_P^{(H_1)}(s)^{q_2}Z_P^{(H_2)}(s)^{q_1}},\label{q1q2}
\end{align}
where $H_1$ and $H_2$ are the subgroups of $G$ whose orders are $q_1$ and $q_2$ respectively.
Then, putting 
\begin{align*}
f(s):=\frac{\zeta_{P_{q_1q_2}}(s)^{q_2}}{\zeta_{P_{q_1q_2}}(q_2s)},
\quad g(s):=\text{RHS of \eqref{q1q2}},
\end{align*}
we have 
\begin{align*}
\frac{f(s)^{q_1}}{f(q_1s)}=g(s).
\end{align*}
from \eqref{q1q2}.
Hence we can obtain the results similar to Theorem \ref{thm1} and \ref{thm2} recursively. 

\end{rem}

\section{Examples}
\subsection{Dedekind zeta functions} 
Let $k$ be an algebraic number field over $\bQ$ such that $[k:\bQ]<\infty$, 
and $K$ a finite Galois extension of $k$. 
When we put $P$ the set of prime ideals of $k$ unramified in $K$ 
and $N$ the norm of $\p\in P$ in $k$, we see that  
$\zeta_P(s)$ is essentially the Dedekind zeta function 
\begin{align*}
\zeta_P(s)=&\prod_{\begin{subarray}{c}\p\in\Prim{(k)}\\ \p\nmid D\end{subarray}}(1-N_k(\p)^{-s})^{-1}
=\zeta_k(s)\prod_{\begin{subarray}{c}\p\in\Prim{(k)}\\ \p\mid D\end{subarray}}(1-N_k(\p)^{-s}),
\end{align*}
where $\Prim(k)$ the set of the prime ideals of $k$ and $D$ is the relative discriminant of $K$ over $k$.
We also put $G:=\Gal(K/k)$ and $\vphi$ the Frobenius automorphism. 
It is well-known that each $L$-function satisfies the properties (i') and (ii). 
Due to the Artin factorization formula, we have
\begin{align*}
Z_P(s)=&\prod_{\begin{subarray}{c}\p\in\Prim{(K)}\\ \p\nmid D\end{subarray}}(1-N_K(\p)^{-s})^{-1}
=\zeta_K(s)\prod_{\begin{subarray}{c}\p\in\Prim{(K)}\\ \p\mid D\end{subarray}}(1-N_K(\p)^{-s}).
\end{align*}
For such zeta functions, we obtain the following results.
\begin{claim}\label{dedekind}
Assume that $G$ is a cyclic group of prime order $q$.
Then the partial zeta function $\zeta_{P_q}(s)$ is analytically continued to $\{\Re{s}>0\}$.
Furthermore, when both $k$ and $K$ are abelian extensions of $\bQ$, 
$\zeta_{P_q}(s)$ has a natural boundary on $\Re{s}=0$.
\end{claim}

The analytic continuation of $\zeta_{P_q}(s)$ in $\{\Re{s}>0\}$ is easily obtained by Theorem \ref{thm1}. 
For proving $\zeta_{P_q}(s)$ has a natural boundary, we prepare the following lemmas.

\begin{lem}\label{distinct} 
Let $\chi^{(1)}_i$ and $\chi^{(2)}_j$ ($1\leq i,j\leq m$) 
be Dirichlet characters respectively modulo $q^{(1)}_i$ and $q^{(2)}_j$ 
such that $\chi^{(1)}_i\neq\chi^{(2)}_j$. Denote by 
\begin{align*}
L_1(s)=\prod_{i=1}^m L(s,\chi^{(1)}_i),\quad L_2(s)=\prod_{i=1}^m L(s,\chi^{(2)}_i).
\end{align*}
Then, for sufficiently large $T>0$, we have 
\begin{align*}
\sum_{\begin{subarray}{c} 0<\Re{\sigma}<1 \\ 0<\Im{\sigma}<T \\ L_1(\sigma)L_2(\sigma)=0 \end{subarray}}
|m_1(\sigma)-m_2(\sigma)|>C\frac{T}{2\pi}\log{T},
\end{align*}
where $C>0$ is a constant and $m_j(\sigma)$ is the multiplicity of $\sigma$ as a zero of $L_j(s)$.
\end{lem}
\begin{proof}
The lemma above for $m=1$ has been proved in \cite{Fu}. 
By applying Bombieri-Perelli's result \cite{BP}, 
we can easily prove that Lemma \ref{distinct} holds for general $m\geq1$.
\end{proof}

\begin{lem}\label{mult}(\cite{Co1}, \cite{Co2} and \cite{Ba}) 
Let $L(s)$ be a Dirichlet $L$-function.
Then, for $M\geq1$ and sufficient large $T>0$, we have
\begin{align*}
\sum_{\begin{subarray}{c}\Re{\sigma}=1/2 \\ 0<\Im{\sigma}<T \\ L(\sigma)=0 \\ m(\sigma)\leq M\end{subarray}}
m(\sigma)>C_M\frac{T}{2\pi}\log{T},
\end{align*}
where $C_M>0$ is a constant which satisfies $C_M=1-O(M^{-2})$ as $M\tinf$.
\end{lem}

\begin{lem}\label{density}(\cite{Mo}) 
Let $L(s)$ be a Dirichlet $L$-function. 
Then, for $0<\alpha<1/2$, we have
\begin{align*}
\sum_{\begin{subarray}{c} 0<\Re{\sigma}<\alpha \\ 0<\Im{\sigma}<T \\ L(\sigma)=0 \end{subarray}}m(\sigma)
\ll T^{5/2\alpha+\epsilon}.
\end{align*}
\end{lem}

\noi{\bf Proof of Claim \ref{dedekind}.}
Let $n:=[K:\bQ]$. 
When $k$ and $K$ are abelian extensions of $\bQ$, the Dedekind zeta functions $\zeta_k(s)$ and $\zeta_K(s)$ 
are expressed by products of Dirichlet $L$-functions. 
Then $g(s)$ is written as 
\begin{align*}
g(s)=\bigg(\disp\prod_{i=1}^n L(s,\chi_i^{(1)})\bigg)\bigg/\bigg(\prod_{i=1}^{n} L(s,\chi_i^{(2)})\bigg).
\end{align*}
From Lemma \ref{distinct}, we see that there exists a constant $C>0$ such that
\begin{align*}
\sum_{\begin{subarray}{c} 0<\Re{\sigma}<1 \\ 0<\Im{\sigma}<T \\ g(\sigma)=0 \end{subarray}}
m(\sigma)>C\frac{T}{2\pi}\log{T}.
\end{align*}
Due to Lemma \ref{mult}, we have
\begin{align*}
\sum_{\begin{subarray}{c} \Re{\sigma}=1/2 \\ 0<\Im{\sigma}<T \\ g(\sigma)=0 \\ m(\sigma)\leq M\end{subarray}}
m(\sigma)
>\tilde{C}_M\frac{T}{2\pi}\log{T},
\end{align*}
where $\tilde{C}_M=m-O(M^{-2})$ as $M\tinf$. 
Then, taking $M$ such that $C+C_M>m$, we obtain 
\begin{align*}
I(T):=\sum_{\begin{subarray}{c} \Re{\sigma}=1/2 \\ 0<\Im{\sigma}<T \\ g(\sigma)=0\end{subarray}}1
>\frac{C+\tilde{C}_M-m}{M}\frac{T}{2\pi}\log{T}.
\end{align*}
Furthermore, from Lemma \ref{density}, we have
\begin{align*}
J_{\alpha}(T):=\sum_{\begin{subarray}{c} 0<\Re{\sigma}<\alpha \\ 0<\Im{\sigma}<T 
\\ g(\sigma)^{-1}=0 \end{subarray}}m(\sigma) \ll T^{5/2\alpha+\epsilon}.
\end{align*}
Then we have
\begin{align*}
\Omega_q(T):=&\#\{\sigma\in\Omega_q\divset 0<\Im\sigma<T\}
>I(T)-\sum_{l\geq1}J_{q^{-l-1}}(q^lT)=O(T\log{T}).
\end{align*}
This implies that the conditions in Theorem \ref{thm2} are satisfied 
and, therefore, the partial zeta function has a natural boundary on $\Re{s}=0$. 
\qed

\subsection{Selberg (Ruelle) zeta functions}
Let $\bH$ be the upper half plane and $\Gamma$ a discrete subgroup of $SL_2(\bR)$ 
such that the volume of $X_{\Gamma}=\Gamma\bsla \bH$ is finite. 
Put $P=\Prim(\Gamma)$ the set of the primitive conjugacy classes of $\Gamma$ and 
$N(p)$ the square of the larger eigenvalue of $p\in\Prim(\Gamma)$. 
Then it is known that $\zeta_P(s)$ coincides the Selberg (Ruelle) zeta function 
\begin{align*}
\zeta_{\Gamma}(s):=\prod_{p\in\Prim(\Gamma)}(1-N(p)^{-s})^{-1}\quad\Re{s}>1
\end{align*}
which satisfies the conditions (i) and (ii) (see, e.g. \cite{He}). 

Fix $\Gamma'$ a normal subgroup of $\Gamma$ of finite index. 
Put $G=\Gamma/\Gamma'$ and $\vphi$ a natural projection from $\Conj(\Gamma)$ to $\Conj(G)$. 
It is known that $L$-functions satisfies (i') and (ii). 
According to \cite{VZ}, we have  
\begin{align*}
\zeta_{\Gamma'}(s):=\prod_{p\in\Prim(\Gamma')}(1-N(p)^{-s})^{-1}\quad\Re{s}>1.
\end{align*}
For the Selberg zeta function, we obtain the following result.

\begin{claim}\label{selberg}
Assume that $G$ is a cyclic group of odd prime order $q$.
Then $\zeta_{P_q}(s)$ is analytically continued to $\{\Re{s}>0\}$.
Furthermore, when both $\Gamma$ and $\Gamma'$ are congruence subgroups of $SL_2(\bZ)$, 
the partial zeta function $\zeta_{P_q}(s)$ has a natural boundary on $\Re{s}=0$.
\end{claim}

\begin{proof}
Similar to the previous section, the analytic continuation of $\zeta_{P_q}(s)$ 
is obtained by Theorem \ref{thm1}.

When $\Gamma$ is a congruence subgroup of $SL_2(\bZ)$, 
due to the determinant expression of the Selberg zeta function (see \cite{Hu} and \cite{Ko}),
we have 
\begin{align*}
\Lambda=\{1/2+ir_j\divset \lambda_j>1/4\}
\cup\bigcup_{l=1}^{h_{\Gamma}}\{\sigma\divset L(2\sigma,\chi_l)=0,0<\Re{\sigma}<1/2,\Im(\sigma)>0\},
\end{align*}
where $\lambda_j=1/4+r_j^2$ is the $j$-th eigenvalue of the Laplacian on $X_{\Gamma}$ with 
the multiplicity $m_j$, 
$h_{\Gamma}$ is the number of cusps for $\Gamma$ 
and $\chi_1,\cdots,\chi_{h_{\Gamma}}$ are Dirichlet characters determined by $\Gamma$. 
Then, similar to the case of the previous section, 
we can easily check that the conditions in Theorem \ref{thm2} are satisfied for $q>2$.
Therefore we obtain Claim \ref{selberg}.
\end{proof}

\begin{rem}
When $X_{\Gamma}$ and $X_{\Gamma}'$ are compact Riemann surfaces, 
it is known that $\Lambda=\{1/2+ ir_j\}$
\begin{align*}
\sum_{\begin{subarray}{c}|r_j|<T \end{subarray}}m_j\sim&
\frac{\vol(X_{\Gamma})}{4\pi}T^2,\qquad \quad m_j\ll \frac{r_j}{\log{r_j}}.
\end{align*}
If it would be known that there exists a constant $\delta>0$ such that 
\begin{align}\label{distinctselberg}
\sum_{\begin{subarray}{c}\sigma=1/2+i\beta \\ 0<\beta<T \\ \{f_1(\sigma)f_2(\sigma)\}^{-1}=0\end{subarray}}
|m_1(\sigma)-m_2(\sigma)|\gg T^{1+\delta}
\end{align}
for $f_1(s):=\zeta_{\Gamma}(s)^q$, $f_2(s):=\zeta_{\Gamma'}(s)$ and 
the orders $m_1(\sigma), m_2(\sigma)$ respectively of the singular points of $f_1(s), f_2(s)$, 
then we could prove that the partial zeta function has a natural boundary on $\Re{s}=0$.
However, no results such like Lemma \ref{distinct} or \eqref{distinctselberg} 
have been obtained for the Selberg zeta functions. 
Hence we cannot presently conclude that whether the partial Selberg zeta function for the compact case 
has a natural boundary on $\Re{s}=0$.
\end{rem}

\subsection{Ihara zeta functions}
Let $X$ be a finite connected $(q+1)$-regular graph with $n$ vertices 
and $Y$ a  finite connected $(q+1)$-regular unramified covering graph of $X$. 
Denote $\Gamma_X$ and $\Gamma_Y$ by the fundamental groups of $X$ and $Y$ respectively, 
and assume that $\Gamma_Y$ is a normal subgroup of $\Gamma_X$.
Put $P$ the set of the equivalence classes of primitive closed backtrackless tail-less cycles.
Then the Ihara zeta function of $X$ are defined by 
\begin{align*}
\zeta_X(u):=\prod_{p\in P}(1-u^{\nu(p)})\quad |u|<1,
\end{align*}
where $\nu(p)$ is the number of the edges in the cycle $p$.
It is known that $\zeta_X(s)$ has the following determinant expression (see, e.g. \cite{Ih}).
\begin{align}
\zeta_X(u)^{-1}=(1-u^2)^{n(q-1)/2}\det(\Id-Au+qu^2\Id),\label{det}
\end{align}
where $A$ is the adjancency matrix of $X$.
From the formula \eqref{det}, we see that $\zeta_X(u)$ has a finite number 
of poles in $\{q^{-1}\leq|u|\leq1\}$. 
Note that the pole at $u=q^{-1}$ is simple. 
For $N(p):=q^{\nu(p)}$, we have $\zeta_P(s)=\zeta_X(q^{-s})$. 
According to \eqref{det}, we see that $\zeta_P(s)$ satisfies the property (i) and (ii). 

Take $G=\Gamma_X/\Gamma_Y$.
It is known that the $L$-functions satisty the conditions (i'), (ii) 
and $Z_P^{(G)}(s)=\zeta_Y(q^{-s})$ (see \cite{ST}). 
Due to \eqref{det}, we see that the singular points of $\zeta_P(s)$ and $Z_P(s)$ 
are periodically distributed. 
Then we can easily prove the following results by using Theorem \ref{thm1} and \ref{thm2}.

\begin{claim}
Assume that $G$ is a cyclic group of odd prime order $q$.
Then $\zeta_{P_q}(s)$ is analytically continued to $\{\Re{s}>0\}$ 
has a natural boundary on $\Re{s}=0$. \qed
\end{claim}

\noindent{\bf Acknowledgement.} 
The author would like to thank to Prof. Nobushige Kurokawa for his advise and helpful comments.

\noindent 
\text{HASHIMOTO, Yasufumi}\\ 
Graduate School of Mathematics, Kyushu University.\\  
6-10-1, Hakozaki, Fukuoka, 812-8581 JAPAN.\\ 
\text{hasimoto@math.kyushu-u.ac.jp}

\end{document}